\documentclass{amsart}
\usepackage{amssymb}
\usepackage{amscd}
\input xypic
\xyoption{all}

\newtheorem{theorem}{Theorem}[section]

\newtheorem{proposition}[theorem]{Proposition}
\newtheorem{corollary}[theorem]{Corollary}
\theoremstyle{definition}
\newtheorem{definition}[theorem]{Definition}

\numberwithin{equation}{section}
\newcommand{\PP}{\mathbb{P}}
\newcommand{\C}{\mathbb{C}}
\newcommand{\pair}[1]{\langle #1 \rangle}
\newcommand{\res}{\operatorname{Res}}
\newcommand{\ind}{\operatorname{ind}}
\newcommand{\HH}{\mathcal{H}}
\newcommand{\alog}{A_\textup{log}}
\newcommand{\aloga}{A_{\textup{log},1}}
\newcommand{\Q}{\mathbb{Q}}
\newcommand{\hdr}{H_{\textup{dR}}}
\newcommand{\sF}{\mathcal{F}}
\begin{document}
\title{A short proof of de Shalit's cup product formula}
\author{Amnon Besser}
\address{Department of Mathematics\\
Ben-Gurion University of the Negev\\
P.O.B. 653\\
Be'er-Sheva 84105\\
Israel
}
\email{bessera@math.bgu.ac.il}

\begin{abstract}
We give a short proof of a formula of de Shalit, expressing the cup
product of two vector valued one forms of the second kind on a Mumford
curve in terms of Coleman integrals and residues. The proof uses the
notion of double indices on curves and their reciprocity laws.
\end{abstract}

\maketitle

\section{Introduction}
\label{sec:intro}

In \cite{deS88} de Shalit proved a formula for the cup product of two
vector valued differential forms on a Mumford curve. This is based on
an earlier partial result of his~\cite{deS89} for two holomorphic
differentials. This formula was later reproved by Iovita and
Spiess~\cite{Iov-Spi01}. The goal of this short note is to give an
alternative short 
proof of de Shalit's formula, based on the theory of the double
index~\cite[Section~4]{Bes98b}.

Let us state de Shalit's result. Let $K$ be a finite extension of
$\Q_p$. Consider a Mumford curve $\HH/\Gamma$,
where  $\Gamma\subset \operatorname{PGL}_2(K)$ is a Schottky group and
$\HH \subset \PP_K^1$ is the rigid analytic space  obtained by
removing the limit points of $\Gamma$.
Let $V$ be a finite dimensional $K$-vector space with a representation
of $\Gamma$. The group $\Gamma$ acts on the space of $V$-valued
differential forms on $\HH$, $\Omega^1(\HH,V) $, by the rule
\begin{equation*}
  \gamma(\sum \omega_i v_i) = \sum (\gamma^{-1})^\ast \omega_i \gamma(v_i) 
\end{equation*}
(compare \cite[1.1]{deS88}). We let it act by the same formula on
spaces of functions.
A $V$-valued differential one-form $\omega$ on $\HH$ with values
in $V$ is $\Gamma$-invariant if $\gamma(\omega)= \omega $ for every
$\gamma\in \Gamma$. It is of the second kind if its
residues (with values in $V$, computed coordinatewise, in any basis),
are $0$ at any point $z\in \HH$. Let $\pair{~}$ be a
$\Gamma$-invariant bilinear form on $V$. The cup product of two
$\Gamma$-invariant $V$-valued one forms of the second kind $\omega$
and $\eta$ can be described by the formula
\begin{equation*}
  \omega\cup  \eta= \sum_{z\in \Gamma\backslash\HH}
  \res_z\pair{F_\omega,\eta}\;,
\end{equation*}
where $F_\omega$ is any primitive of $\omega$ locally near $z$, which
exists (formally) because of the residue of $\omega$ at $z$ is $0$, and
is independent of the choice of the primitive because the residue of
$\eta$ at $z$ is $0$. Note that the expression to be summed indeed
depends only on $z$ modulo $\Gamma$.

An open annulus is a rigid space isomorphic to the space $s<|z|<r
$. An orientation on an annulus may be described as a choice of a
parameter $z$ as above, with two parameters considered equivalent if
they give the residue, as defined below. An annulus together with an
orientation is called an oriented annulus. A differential form
$\omega$ on an oriented annulus $e$ has a residue $\res_e \omega $
such that $\res \sum a_i z^i dz= a_{-1} $. It can be shown that there
are only two orientations, giving residues differing by multiplication
by $-1$. By choosing a basis for $V$
the residue extends easily to $V$ valued differential forms.

We now recall~\cite[Definition~2.5]{deS89} that the action of $\Gamma$
on $\HH$ has a good
fundamental domain in the following sense: There are pairwise disjoint
closed $K$-rational discs $B_i$ and $C_i$, $i=1,\ldots, g $ and open
annuli $b_i$, $c_i$, and elements $\gamma_i \in \Gamma$, such that the
following holds:
\begin{enumerate}
\item The $\gamma_i$ freely generate $\Gamma$.
\item The unions $B_i\cup b_i$ and $C_i\cup c_i$ are open discs, still
  pairwise disjoint.
\item For each $i$, $\gamma_i$ maps $B_i$ isomorphically onto the
  complement of $C_i\cup c_i$ and $b_i$ isomorphically onto $c_i$.
\item The complement of $\bigcup_i (B_i\cup b_i \cup C_i) $ is a
  fundamental domain for $\Gamma$.
\end{enumerate}
We give the annuli $c_i$ and $b_i$ the orientation given by the discs
$C_i$ and $B_i$ respectively, i.e., one given by parameters extending
to $C_i \cup c_i$ and taking the value $0$ on $C_i$ (respectively with
$b_i$ and $B_i$). Thus, $c_i$ is oriented in the same way
as in~\cite[1.5]{deS88} while $b_i$ is oriented in the reversed direction
to loc.\ cit. (the $b_i$'s do not show up in the formula). With this
choice, $\gamma_i:b_i \to c_i $ is orientation reversing.

de Shalit's formula involves Coleman integration of holomorphic
$V$-valued $1$ forms. While this can be described in a completely
elementary way since we are dealing with subdomains of the projective
line~\cite[P. 41]{Ger-Put80}, we will use the more involved theory of
Coleman~\cite{Col-de88} 
and adapt it to our case by choosing a basis of $V$ and then integrate
coordinate by coordinate. This is clearly independent of the choice of
a basis because Coleman integration is linear (up to constant). The
key property of Coleman integration is its functoriality. It
immediately implies that from the property $\gamma \omega= \omega $ we
may deduce that for any $\gamma\in \Gamma$ the function $\gamma(F_\omega)-
F_\omega$ is constant. We can now state the main theorem.
\begin{theorem}[{\cite[Theorem~1.6]{deS88}}]\label{main}
  With the data above we have
  \begin{equation*}
    \omega \cup \eta = \sum_i  \pair{\gamma_i  F_\omega- F_\omega,
  \res_{c_i} \eta} - \pair{\res_{c_i} \omega,  \gamma_i F_\eta - F_\eta }\;.
  \end{equation*}
\end{theorem}
The main ingredient in the present proof is the theory of double
indices and their reciprocity laws on
curves~\cite[Section~4]{Bes98b}. We need a
very easy extension of this theory to vector valued differential
forms. Once this has been described, the proof is an easy computation.

We would like to thank Michael Spiess for suggesting this project. The
author is supported by a grant from the Israel Science Foundation.

\section{Double indices of vector valued differential forms}
\label{sec:double}

In this section we describe a rather straightforward generalization of
the theory of double indices~\cite[Section~4]{Bes98b} to the case of
vector valued one
forms. The extension is fairly trivial since we consider only constant
coefficients. We work over $\C_p$ for convenience.

Let $A$ be either the field of meromorphic functions in the variable
$z$ over $\C_p$ or the
ring of rigid analytic functions on an annulus  $\{r<|z|<s\}$ 
over $\C_p$. Let $\alog:=A[\log(z)]$ and let $\aloga\subset \alog$ be
the subspace of $F\in \alog $ which are linear in $\log(z)$, a
condition which is
equivalent to $dF\in A dz$.
\begin{definition}{\cite[Proposition~4.5]{Bes98b}}
  The double index, $\ind(~):\aloga \times \aloga \to \C_p $ is the
  unique antisymmetric bilinear pairing such that $\ind(F,G) =\res F
  dG$, whenever $F\in A$.
\end{definition}
Suppose now that $C$ is a proper smooth curve over $\C_p$ with good
reduction, and that $U$ is a rigid analytic space obtained from $C$ by
removing discs $D_i$ of the form $|z_i|\le r $, with $r<1$,  where the
reduction of $z_i$ is a
local parameter near a point $x_i$ of the reduction. Let us call these
domains \emph{simple domains}. To the disc $D_i$
corresponds the annulus $e_i$ given by the equation $r< |z_i| < 1$,
which is contained in $U$ and oriented by $z_i$.

Choose a branch of the $p$-adic logarithm. Given a rigid one form
$\omega\in \Omega^1(U) $, Coleman's theory 
provides us with a unique up to constant, locally analytic function
$F_\omega$ on $U$ with the property that $d F_\omega = \omega
$. Restricted to the annuli $e_i$ these clearly belong to $\aloga$ and
one can therefore define, for two such functions $F_\omega $ and
$F_\eta$ the double index $\ind_{e_i}(F_\omega,F_\eta) $. It follows
from~\cite[Lemma~4.6]{Bes98b} that this index depends only on the
orientation. One of the
main technical results of~\cite{Bes98b} is the following.
\begin{proposition}[{\cite[Proposition~4.10]{Bes98b}}]
  We have $\sum_i \ind_{e_i}(F_\omega,F_\eta) = \Psi(\omega)\cup
  \Psi(\eta) $, where $\Psi: \hdr^1(U) \to \hdr^1(C) $ is a certain
  projection.
\end{proposition}
We will only need the following immediate Corollary, which follows
because  $\hdr^1(\PP^1/\C_p)=0 $.
\begin{corollary}\label{p1cor}
  Suppose that $C=\PP^1 $. Then, in the situation above, $\sum_i
  \ind_{e_i}(F_\omega,F_\eta)=0$.
\end{corollary}
We can now extend the theory to vector valued differential forms in a
rather trivial way. Suppose we are given a finite dimensional
$\C_p$-vector space with a bilinear form $\pair{,} $.
\begin{definition}
  Chose bases $\{v_i\}$ and $\{u_i\}$ for $V$. Suppose that the
  $V$-valued Coleman
  functions $F_\omega$ and $F_\eta$ are written as
  \begin{align*}
    F_\omega &= \sum F_{\omega_i} v_i\;,\\
    F_\eta &= \sum F_{\eta_i} v_i\;.
  \end{align*}
  Then, the local index $\ind_e(F_\omega,F_\eta) $ is given by
  \begin{equation*}
    \ind_e(F_\omega,F_\eta) = \sum_{i,j}
    \ind_e(F_{\omega_i},F_{\eta_j}) \pair{v_i,u_j}\;.
  \end{equation*}
\end{definition}
It is easy to check that this definition does not depend on the choice
of bases. An easy consequence of the definitions is the following.
\begin{proposition}\label{antivec}
  Suppose that $\res_e \omega=0$. Then $\ind_e(F_\omega,F_\eta)=
  \res_e\pair{F_\omega,\eta} $
  while $\ind_e(F_\eta,F_\omega)=-\res_e\pair{\eta,F_\omega} $.
\end{proposition}

We now restrict to the case $C=\PP^1$ but consider more general
subdomains $U$,  obtained by removing closed discs $D_i =
|z-\alpha_i|=r_i$, including the case of removing a point when $r_i=0$. For
each $i$ we consider an annulus $e_i$ in $U$ surrounding $D_i$, in
such a way that the open discs $D_i \cup e_i$ are still disjoint We
will call the $e_i$ the \emph{annuli ends} of $U$. It
is easy to see that $U$ can be obtained by gluing simple domains $U' \in
\PP^1$ along 
annuli. Note that the $U'$'s are glued along annuli with
reversed orientations. Given $\omega\in \Omega^1(U,V)$ one can define
its Coleman 
integral $F_\omega$ first on each of the $U'$'s as before and then by
adjusting constants along the annuli. The intersection graph of the
$U'$'s is a tree so there is always a way of choosing an integral
globally. This construction coincides with the definition of Coleman
integrals in~\cite{Ger-Put80}.
\begin{proposition}\label{easyrecip}
  In the situation described above we have, for any rigid $V$-valued
  one-form on $U$,   $\sum_i \ind_{e_i}(F_\omega,F_\eta) = 0$.
\end{proposition}
\begin{proof}
The case $V$  trivial and $U$ simple is Corollary~\ref{p1cor}. We next
consider the
case $U=U'_1 \cup U'_2 $ with $U'_1$ and $U'_2$ glued along an annulus
$e$. Since $e$ has reversed orientations when considered in $U'_1$ and
$U'_2$, the double index $\ind_e $ has a reverse sign in these two
cases by~\cite[Lemma~4.6]{Bes98b}. Thus, the result for $U$ follows
from those for $U'_1$ and $U'_2$. Now, the case of a general
$U$, still with trivial $V$, follows immediately.  The general case
follows by choosing bases.
\end{proof}
\begin{proposition}\label{2.5}
  Let $e$ be an annulus in $\HH$ and let $\gamma\in \Gamma$. For
  $\omega\in \Omega^1(e,V) $ let $F_\omega $ be its integral. Then
  $\gamma F_\omega$ is a Coleman integral of $\gamma(\omega) $ on
  $\gamma(e)$, furthermore, If $\eta$ is another such form, then we
  have
  \begin{equation*}
    \ind_e(F_\omega,F_\eta) = \pm \ind_{\gamma(e)}
    (\gamma(F_\omega),\gamma(F_\eta))\;,
  \end{equation*}
  depending on whether $\gamma$ is orientation reversing or saving.
\end{proposition}
\begin{proof}
We choose a basis $\{v_i\}$ of $V$ and we let $\{u_i\}$ be the dual
basis with respect to $\pair{,}$. Then, since $\pair{,} $ is
$\Gamma$-invariant, the bases
$\{\gamma(v_i)\}$ and $\{\gamma(u_i)\}$ are also dual to each
other. This implies that if $F_\omega= \sum f_i v_i$ while $F_\eta=
\sum g_i u_i$, then 
\begin{align*}
  \ind_e(F_\omega,F_\eta)&= \sum_i \ind_e(f_i,g_i)\\
 \intertext{and}
  \ind_{\gamma(e)} (\gamma(F_\omega),\gamma(F_\eta)) &= \sum_i
  \ind_{\gamma(e)}((\gamma^{-1})^\ast f_i,(\gamma^{-1})^\ast g_i)\;.
\end{align*}
But by \cite[Lemma~4.6]{Bes98b} we have, for each $i$,
\begin{equation*}
  \ind_e(f_i,g_i)= \pm \ind_{\gamma(e)}((\gamma^{-1})^\ast
  f_i,(\gamma^{-1})^\ast g_i)\;,
\end{equation*}
depending on whether $\gamma^{-1}$ is orientation reversing or
preserving, and the result follows immediately from this.
\end{proof}

\section{The proof}
\label{sec:proof}

\begin{proof}[Proof of Theorem~\ref{main}]
By the remark following Equation (5) in~\cite{deS88} we may assume
that $b_i$ and $c_i$ contain no poles of $\omega$ and $\eta$.
Consider the domain $\sF = \PP^1- \bigcup_i (B_i \cup C_i) $, which is
of the type considered in Section~\ref{sec:double}, and its annuli
ends are the $b_i$ and $c_i$. It follows from the description of the
fundamental domain for $\Gamma$ that $\sF - \bigcup_i (c_i\cup b_i) $ contains exactly one
out of every $\Gamma$ class of every singularity of either forms. It
follows that
\begin{equation*}
  \omega \cup \eta = \sum_{x\in \sF} \res_x \pair{F_\omega,\eta} =
  \sum_{x\in \sF}
  \ind_{x}(F_\omega,F_\eta) = - \sum_i (\ind_{b_i}(F_\omega,F_\eta)+
  \ind_{c_i}(F_\omega,F_\eta))
\end{equation*}
where the last equality follows from Proposition~\ref{easyrecip}. We now
observe that since $\gamma_i$ is orientation reversing we have by
Proposition~\ref{2.5} that 
$\ind_{b_i}(F_\omega,F_\eta)=-\ind_{c_i}(\gamma_i F_\omega,\gamma_i
F_\eta) $. Therefore
\begin{align*}
  &\phantom{=}-(\ind_{b_i}(F_\omega,F_\eta)+  \ind_{c_i}(F_\omega,F_\eta))\\&=
  \ind_{c_i} (\gamma_i F_\omega, \gamma_i  F_\eta)-\ind_{c_i}(F_\omega,F_\eta)
\\ &=
  \ind_{c_i} (\gamma_i F_\omega-F_\omega, \gamma_i F_\eta)+
  \ind_{c_i}(F_\omega,\gamma_i F_\eta-F_\eta)
\\ &= \res_{c_i} \pair{\gamma_i F_\omega-F_\omega, \gamma_i \eta}-
  \res_{c_i}\pair{\omega,\gamma_i F_\eta-F_\eta}\quad \text{by
  Proposition~\ref{antivec}}
\\ &= \pair{\gamma_i F_\omega-F_\omega, \res_{c_i}  \eta}-
  \pair{\res_{c_i}\omega,\gamma_i F_\eta-F_\eta}\;.
\end{align*}
The theorem follows immediately.
\end{proof}

\end{document}